\documentclass[11pt]{article}
\usepackage{amsmath}
\usepackage{amssymb}
\usepackage{hyperref}
\usepackage{caption}
\usepackage{amsthm}
\usepackage{comment}

\newtheorem*{GI}{Gronwal Inequality}
\newtheorem{theorem}{Theorem}
\newtheorem*{proof1}{An Alternative Proof for Theorem 1}
\newtheorem*{proof2}{Proof of theorem 2}
\newtheorem{rem}{Remark}
\newtheorem*{question}{Question}
\newtheorem*{example}{Example}

\begin{document}

\title{A Dynamical Approach to Quasi-Analytic type Problems }
\author{Ali Taghavi \\ Department of Mathematics, Qom university of Technology\\
Qom, Iran,\\
 taghavi@qut.ac.ir}

\maketitle
\begin{abstract}
In this paper we give an alternative proof for a vanishing result about flat functions proved in G.Stoica, "When must a flat function be identically zero", The American Mathematical Monthly 125(7)648-649, 2018. With a dynamical approach we give  a generalization of this result to multidimensional variables.
\end{abstract}

\section*{Introduction}
Let $M=(m_k)$ be a sequence of positive real numbers. We consider the class $C(\{m_n\})$   of  all smooth functions on an interval $[a,b]$ such that for all $x\in [a,b]$  we have $|f^{(n)}(x)|\leq b^n.m_n$ where $b$ is a constant which depends on $f$. A  class $C(\{m_n\})$ of smooth functions is called a quasi analytic class if for every $f\in C(\{m_n\})$ if $f^{(n)}(a)=0,\quad n=0,1,2,\ldots $ at some point $a$ then $f$ must be identically $0$ on $[a,b]$. This concept generates several results including Denjoy-Carlman theorem. See \cite{Car,Den}. This theory was a motivation for \cite{stoica} to consider the following question: Under what functional differential inquality imposed on a flat function $f$ and its derivative  can one conclude that $f$ identically vanishe?. So this is an strong motivation to extend the concept of quasi analytic functions via consideration of various functional differential inequalities.  In this direction  George Stoica   proved  the following interesting result in  \cite{stoica}:\newpage

\begin{theorem} \label{GS}  Let f be a $C^1$ real-valued function on $[0, 1]$, infinitely differentiable at
$x = 0$, and such that $f^{(n)}(0) = 0,\;\forall n\in \mathbb{N}\cup \{0\}$. If
$|xf'(x)| \leq C|f (x)|$ for some $C > 0$ and every $x \in [0, 1]$, then $f (x) = 0$ for every $x \in  [0,1]$\end{theorem}

In \cite{HLOU} the author gives a simplified proof for the above theorem.

In this paper, apart from giving an alternative proof for this theorem,  we observe that  there is a dynamical nature for the formulation of this  theorem. The method of proof of the following theeorem 2, which generalizes \autoref{GS} to the  multdimensional variables, represents  this dynamical feature.

\begin{theorem} \label{TT}   Assume that $U$ is a disc around the origin $0\in \mathbb{R}^n$. Let $h:U \to \mathbb{R}^n$ with   $h=(h_1,h_2,\ldots,h_n)$ be a $C^1$ map with the inner product condition $\langle h(x), x \rangle >0, \forall x \in U\setminus \{0\}$.   We further assume that all eigenvalues of the Jacobian matrix $Jh(0)$ have  positive  real part.
Assume that $f:U \to \mathbb{R}$ be a $C^1$ function which is flat at $0$. If we have $|\sum_{i=1}^n h_i\partial f/\partial x_i| \leq c |f(x|$ for some constant $c$ then $f$ is identically $0$ on $U$.\end{theorem}

\begin{example}
Put $x=(x_1,x_2,\ldots x_n)\in \mathbb{R}^n$. Let $f:\mathbb{R}^n \to \mathbb{R}$ be a smooth function which is  flat  at $0$. Assume  that $|\sum_{i=1}^{n} x_i \frac{\partial f}{\partial x_i}|\leq c|x|$. Then $f(x)=0,\forall x\in \mathbb{R}^n$. To prove this fact we apply \autoref{TT} to $h(x)=x$. The Jacobian matrix $Jh=I_n$ whose only eigenvalue is $1$. Moreover $\langle h(x), x \rangle >0.\quad \forall x\neq 0$
\end{example}

\section*{ Backgrounds from Dynamics and Proofs of the main Theorems}

In this section we first recall  the Gronwal inequality and flat functions. Then we provide necessary definitions from the theory of differential equations and dynamical systems.

\begin{GI}
 Let $\beta$ and $u$ be real-valued continuous functions defined on an interval $[a,b]$. If $u$ is differentiable in $(a,b)$ and satisfies the differential inequality $u'(t)\leq \beta(t)u(t).\;\;\forall t\in [a,b]$ then $u(t)\leq u(a) \exp{\int_a^t \beta(s)ds}$
 for all $t\in [a,b]$.\end{GI}

Let $f$ be a  real valued function defined in a neighborhood of $0\in \mathbb{R}$. We say that $f$ is a flat  function at the origin if $f^{(n)}(0)=0,\;\forall n\in \mathbb{N}\cup\{0\}$. Its multivariable version  is the following:\\
 Assume that $U$ is an open subset of $\mathbb{R}^n$ containing $0$. A $C^1$ function $f:U \to \mathbb{R}$ is called a flat function at the origin, if $f(0)=0$ and  all its partial derivatives at $0$ exist and vanishes. Namely $\frac{\partial ^m f}{\partial x_i^p\partial x_j^q}(0)=0,\;\;\forall m\in \mathbb{N};, 1\leq i,i\leq n\;\; p+q=m$. This is equivalent to say that $f$  is infinitely differentiable at the origin and \begin{equation} \label{equivalence}  \lim_{x\to 0}  \frac{|f(x)|}{|x|^k}=0\quad \forall k\in \mathbb{N}\end{equation}

Now we shall justify the dynamical nature of Theorem 1. First we introduce some dynamical preliminaries  which can be found in  \cite{hirsch-smale} :\\

Let $U$ be an open subset of $\mathbb{R}^n$ and $X:U\to \mathbb{R}^n$ be  a vector field  with $X=(P_1,P_2,\ldots P_n)$. By existence and uniqueness theorem of the theory of ordinary differential equations we get a flow $\phi_t$ associated to solution curves of \begin{equation} \label{ET}  x'=X(x) \end{equation}  namely $\phi_t(x)$ is the solution curve of ~\eqref{ET}  staring at point $x$,  See\cite[Chapter 8]{hirsch-smale}. The maximal interval of solution starting at $x$, denoted by $I(x)$, is the maximal interval around $t=0$ where the solution $\phi_t(x)$ can be defined on $I(x)$. The orbit of a point $x$ is defined as $$\mathcal{O}(x)=\{\phi_t(x)\mid t\in I(x)\}$$ The positive and negative semi orbits $\mathcal{O}^{\pm}$ defined as $$\{\phi_t(x)\mid t\in I(x),\text{and}\quad t\geq 0(t\leq 0,\quad resp.\}$$  A singularity of $X$ is a point $a\in U$ with $X(a)=0$. A singularity $a$  of $X$ is called a hyperbolic sink (source, resp.) if  all eigenvalues of $JX(a)$, the Jacobian matrix of  $X$ at $a$, have negative(positive, resp.) real part. The following theorem about hyperbolic sinks and source plays a crucial role in our paper. Its proof can be found  in \cite{hirsch-smale}:
\begin{theorem}\label{SS}

Let $a\in U$ be a sink of equation $\dot x =X(x)$ with flow $\phi_t$. Then there exist an open set $V\subset U$ containing $a$ and  constant $\theta>0,
;\lambda>0$  such that $|\phi_t(x)-a|\leq \theta e^{-\lambda t}|x-a|,\quad \forall x\in V,\quad \forall t>0$
\end{theorem}

Let $X=(P_1,P_2,\ldots P_n)$ be a vector field on an open set $U \subset \mathbb{R}^n$. Assume that $f:U \to \mathbb{R}$ is a $C^1$ function. We define $X.f=\sum_{i=1}^n P_i\partial f/\partial x_i $. Then $X.f$ is the derivation of $f$ along solution curves of $X$. Namely \newline $X.f(z)=\frac{d}{dt}(f\circ\phi_t(z))|_{t=0}$. The linear  operator $D$ with $D(f)=X.f$ is called "\emph{Derivational operator associated to X}".\\

Assume that a vector field $X$ has  a singularity at a point $a$. Assume that a function $f$ is defined in a neighborhood of $a$ with a local minimum at point $a$. We say that $f$ is  a Liapunov function for $X$ if $X.f<0$ in a deleted neighborhood of $a$. Existence of Liapunov functions around a singularity $a$ implies that all solution curves $\phi_t(x)$ starting points $x$ in a small neighborhoods of $a$ tends to $a$ as $t$ goes to infity. This concept is called stability. See \cite[Chapter 9]{hirsch-smale}.

\begin{rem}

With these notations we provided so far, the \autoref{GS} can be read in the following dynamical language:\\

\texttt{\emph{Let $X$ be the vector field  $X(x)=x$. Assume that $f$ is a flat function at the origin. Then if $|X.f(x)|\leq C|f(x)|$ for some constant $C$ then $f$ vanishes identically.}}\\
So it is natural to consider  a dynamical nature for \autoref{GS}. Such a consideration enabled us to generalize this theorem to higher dimensional case.
 \end{rem}
 Now we shall provide a new proof for the result in \cite{stoica}.
 \begin{proof1}

On the contrary, assume that $f$ is not identically zero. Then there exist a right\_isolated zero for $f$. This means that there exist a point $x_0\in [0, 1)$ with $f(x_0)=0$ such that \begin{equation} \label{H}f(x)\neq 0 \qquad \forall  x\in (x_0,x_0+\delta]\end{equation}  for some $\delta >0$. Without lose of generality we may assume that $f$ is positive on the open interval $(x_0, x_0+\delta]$, otherwise we replace $f$ by $-f$. First assume that $x_0\neq 0$. We define $u(t)=f(x_0e^t)$ for $t\in [0, \epsilon)$ where $\epsilon=\ln \frac{x_0+\delta}{x_0}$. Now  we differentiate $u(t)$ and apply the assumption $|xf'(x)|\leq C|f(x)|$ of Theorem 1. Then we have \newline \begin{equation}\label{EEE}
  \begin{aligned}
   u'(t)&=x_0e^tf'(x_0e^t) \leq |x_0e^tf'(x_0e^t)|\\& \leq C|f(x_0e^t)|=Cf(x_0e^t)=Cu(t)
  \end{aligned}
\end{equation}

 By  Gronwal inequality we get $0\leq u(t)\leq u(0)e^{ct}=0$. This means that $f$ vanishes in a right neighborhood of $x_0$. This  contradicts to the assumption  ~\eqref{H}.

Now assume that $x_0=0$. So $f$ does not vanish on $(0, \delta]$. For every fixed $x\in (0,\delta)$  we define $u(t)=f(xe^t),\;\;t\in [0,\ln \frac{\delta}{x}]$. We apply again the Gronwal inequality to \eqref{EEE} so we obtain $f(xe^t)\leq f(x)e^{Ct},\;\;\forall t\in [0, \frac{\delta}{x}]$. We substitute $t= \frac{\delta}{x}$ in the latter inequality so we have $0<f(\delta)\leq f(x)\left(\frac{\delta}{x}\right)^C\Rightarrow f(x)>\left(\frac{f(\delta)}{\delta}\right )^Cx^C$. This obviously contradicts with flatness of $f$ at the origin. This completes the proof of Theorem 1.

\end{proof1}

Our proof of the above theorem can be generalized to higher dimensions. This generalization is given in the proof of \autoref{TT} as follows:

\begin{proof2}
The method of proof is based on usage of \autoref{SS} and also on the  same method we used in the proof of  \autoref{GS}. We apply the same method to a given typical orbit or semi orbit $\mathcal{O}(p)$ or $\mathcal{O}^{\pm}(p)$ of the vector field $h$ or its negative direction $Y=-h$.  The vector field $h$ has a source at the origin hence the vector field  $Y=-h$ has a sink at $0$. We denote by $\phi_t, \psi_t$ the flow of $Y$ and $h$ respectively. Obviously  $\psi_t(x)=\phi_{-t}(x)$. For every $p\in U, \lim_{t\to +\infty} \phi_t(p)=0$ since $$Y.(x_1^2+x_2^2+\ldots+x_n^2)=-h.(x_1^2+x_2^2+\ldots+x_n^2)=-2\langle h(x).x\rangle<0$$ when $ x=(x_1,x_2,\ldots,x_n)\in U\setminus\{0\}$. So the function $(x_1^2+x_2^2+\ldots+x_n^2)$ is strictly decreasing along the solution curves of $Y$. In the other words the function $(x_1^2+x_2^2+\ldots+x_n^2)$ is a Lyapunov function for the vector field $Y$. To prove \autoref{TT} we assume on the contrary that there exist a point $p\in U$ with $f(p)\neq 0$.  Since $Y$ has a sink at origin which attracts all solution curves  of  disc $U$, the maximal interval of  solution $I(p)$ associated to vector field $h=-Y$ is in the form \begin{equation} \label{ddd} I(x)=(-\infty, \omega(x))\end{equation}. We define a function $u(t): I(x)\to \mathbb{R}$  with $u(t)=f\circ \psi_t(p)$. So $u$ is not identically zero. Furthermore the condition $|\sum_{i=1}^n h_i\partial f/\partial x_i| \leq c |f(x|)$ implies that $u$ satisfies  $u'(t)\leq cu(t)$. So one of the following two situations may occur:\\
\begin{enumerate}
\item \label{GGGG} There exist a right\_isolated root $t_0$ for $u$. This means that there exist a point  $t_0\in I(x)$ with $u(t_0)=0$ but $u$ does not vanish on a right neighborhood $(t_0,t_0+\delta)$\\
\item \label{ff} The function $u$ does not vanish on an interval $(-\infty, m)$ for some  $m<0$
\end{enumerate}

 The first \autoref{GGGG} leads us to contradiction by Gronwal inequality. the function $u(t)$ satisfies $u(t)\leq e^{ct}u(t_0)=0 $ which contradict to non vanishing assumption on $(t_0,t_0+\delta)$.  Note that a similar argument worked in the proof of \autoref{GS}. We consider the second \autoref{ff}. According to \autoref{SS} there exist a neighborhood $V$ around the origin and positive constants $\theta, \lambda$  such that $|\phi_t(x)|\leq \theta e^{\lambda t}|x|,\quad \forall x\in V,\;\forall t>0$. We fix a point $q\in \mathcal{O}^{-}(p)$ which belongs $V$. For every $x=\phi_t(q),\quad t>0$ we shall apply the Gronwal inequality to function $w(s)=f(\psi_s(x))$ defined on $[0,t]$. Again  the assumptions of \autoref{TT} imply that $dw/ds=\dot w\leq cw$ hence $w(s)\leq e^{cs}w(0)$. So we have \begin{equation*} \begin{aligned} f(q)=w(t)\leq e^{ct}w(0)=e^{ct}f(\phi_t(q)) \Rightarrow \\
  f(q)\left (\frac{|\phi_t(q)|}{|q|}\right )^{c/\lambda}\leq f(q){(e^{-\lambda t})}^{c/\lambda}\leq f(q)e^{-ct}\leq f(\phi_t(q)) \end{aligned} \end{equation*}

  \end{proof2}

 So we obtain $$k|\phi_t(p)|^{c/\lambda}\leq f(\phi_t(p))\qquad \forall t>0$$ where the constant $k$ is $k=\frac{f(q)}{|q|^{c/\lambda}}$  Thanks to the  equivalent formulation  of flatness property `\eqref{equivalence}  this contradicts to the fact that  $f$ is flat at origin since $\phi_t(q)$ tends to the hyperbolic sink $0$ as $t$ goes to $-\infty$. This completes the proof of \autoref{SS}

 \section*{Discussions and further Reseaches}

 In this paper we presented a differential operator interpretation for  the assumption $|xf'(x)|\leq c|x|$  in\cite{stoica}. We considered the differential operator $D(f)=xf'(x)$ so we translated  the main problem  as follows: Every flat function $f$ with $|D(f)(x)|\leq c|x|$ must be vanished identically.
 This situation leads us to the following definition and questions:\\

 let $D$ be  a  differential operator   which acts on the space of all smooth functions on $\mathbb{R}^n$ or any arbitrary manifold. We say that $D$ is a $G.S.$ operator at point $p$ if for every flat function locally defined around $p$ which satisfies $|D(f)(x)|\leq c|x|$ then $f$ must  vanish identically in a neighborhood  around $p$. We learn from \cite{stoica} that the differential operator $D(f)=xf'(x)$ is a $G.S.$ operator at origin. Moreover we  prove in this paper that the operator $D(f)=\sum_{i=1}^n P_i(\partial f/\partial x_i) $ is a multidimensional $G.S.$ operator at the origin if $0$ is a hyperbolic sink(or source) for the vector field $X=(P_1,P_2,\ldots,P_n)$. These situations suggest  some   proposals as follows:
 \begin{itemize}
 \item  Introducing more examples of $G.S.$ differential operators and a possible classification of all such kind of operators

 \item Consideration of the previous proposal for the particular operator $\Delta$, the Laplace operator on a Riemannian manifold. Determination of all Riemannian manifold $M$ for which the Laplacian is a $G.S.$ operatore at all points. On the opposite extreme one can think to classification of all Riemannian manifold $M$ for which the operator $\Delta$ is a $G.S.$ operator at no point of $M$

 \item To a vector field $X$ with corresponding  differential operator $D(f)=X.f$, we associate the set  of all points $p$ such that $D$ is a $G.S.$ operator at $p$. So it would be interesting to study this set from a dynamical point of view.\\

     \end{itemize}

     Apart from the above geometric proposals arising from the concept introduced in \cite{stoica}, in the following remark we present a different question with a functional analysis nature:

\begin{rem}

The assumption $|D(f)(x)|\leq c|x|$ is suggesting a concept stronger than usual bounded ness of operator $D$ when we consider the sup norm on an appropriate function space. So this is a motivation to ask the following question when we reduce this strong continuity to the standard boundedness of $D$;\\

Let $V=\{f\in C^{\infty}[0,1]\mid f\quad \text{is flat at 0} \}$. We equip $V$ with the norm $|\quad |_{\infty}$. Is there an infinite dimensional subspace $W\subset V$ which is invariant under $D(f)=xf'(x)$ and $D$ is a  bounded operator on $W$?

\end{rem}

We  observe that the concept of quasi analytic function was the main motivation for consideration of  the differential inequality $|xf'(x)|\leq c|x|$ in \cite{stoica}. The classical concept of quasi analytic functions involves a countable functional inequalities. So a general and natural question is the following:

\begin{question}
Can we study a typical quasi analytic class of functions via a unique functional differential inequality? If there are some partial affirmative answer to this question,  according to existence of  various dynamical interpretations for a functional differential inequality, we ask the next  question:  Is it reasonable and promising  to follow   some dynamical approaches to  classical quasi analytic  problems?

\end{question}

\end{document}